\documentclass[12pt]{article}
\usepackage{epsfig}
\usepackage{latexsym}
\begin{document}

\title{\Large{\bf Circular colorings, orientations, and weighted digraphs}}

\author{Hong-Gwa Yeh\thanks{Partially supported by the
National Science Council of R.O.C. under
grant NSC95-2115-M-008-006.} \\
\normalsize   Department of Mathematics\\
\normalsize   National Central University\\
\normalsize   Jhongli City, Taoyuan 320, Taiwan}

\date{\small December 2006}

\maketitle

\newtheorem{theorem}{Theorem}
\newtheorem{lemma}[theorem]{Lemma}
\newtheorem{corollary}[theorem]{Corollary}
\newtheorem{definition}[theorem]{Definition}
\newtheorem{proposition}[theorem]{Proposition}
\newtheorem{conjecture}[theorem]{Conjecture}
\newcommand{\qed}{\hfill $\Box$ }
\def\G{{\vec{G}}}
\def\E{{\vec{E}}}
\def\C{{\vec{C}}}
\def\1{\mbox{\boldmath $1$}}
\def\T{{\mathcal{T}}}
\def\f{{f_\mathcal{T}}}
\def\ff{{f_{\mathcal{T}'}}}
%
%
\newenvironment{proof}{
\par
\noindent {\bf Proof.}\rm}%
{\mbox{}\hfill\rule{0.5em}{0.809em}\par}
%
%
\newenvironment{proofif}{
\par
\noindent {\bf Proof of the `if' part of Theorem 1.}\rm}%
{\mbox{}\hfill\rule{0.5em}{0.809em}\par}
\newenvironment{proofonlyif}{
\par
\noindent {\bf Proof of the `only if' part of Theorem 1.}\rm}%
{\mbox{}\hfill\rule{0.5em}{0.809em}\par}
\newenvironment{proofcorollary2}{
\par
\noindent {\bf Proof of Corollary 2.}\rm}%
{\mbox{}\hfill\rule{0.5em}{0.809em}\par}
\newenvironment{proofcorollary3}{
\par
\noindent {\bf Proof of Corollary 3.}\rm}%
{\mbox{}\hfill\rule{0.5em}{0.809em}\par}

 \baselineskip=20pt
\parindent=1cm

 %
 %
 %
 %
 \begin{abstract}
 In this paper we prove that if a weighted symmetric digraph $(\G,c)$
 has a mapping $T:E(\G)\rightarrow \{0,1\}$ with $T(xy)+T(yx)=1$
 for all arcs $xy$ in $\G$  such that for each dicycle $C$ satisfying
 $0< |C|_c(\mbox{mod $r$})< \max_{xy\in E(\G)}c(xy)+c(yx)$ we have
 ${|C|_c/ |C|_T}\leq r$,
 then $(\G,c)$ has a circular
 $r$-coloring.
 Our result generalizes the work of Zhu (J.~Comb.~Theory, Ser.~B, 86 (2002) 109-113)
 concerning the $(k,d)$-coloring of a graph, and thus is also a
 generalization of a corresponding result of Tuza
 (J.~Comb.~Theory, Ser.~B, 55 (1992) 236-243). Our result also
 strengthens a result of Goddyn, Tarsi and Zhang (J.~Graph Theory 28 (1998) 155-161)
 concerning the relation between orientation and the
 $(k,d)$-coloring of a graph.
 \end{abstract}

%
%
%
%
\section{Introduction}
 A {\em weighted digraph} is denoted by $(\G,c)$, where $\G$ is a
 digraph, and $c$ is a function which assigns to each arc of $\G$ a
 positive real number. For simplicity of notation, the arc $(u,v)$ is
 written as $uv$, and $c(uv)$ is written as $c_{uv}$.
 If arcs
 $uv$, $vu$ both exist and do not exist for all vertices $u,v$ in
 $\vec{G}$ then
 $(\vec{G},c)$ is said to be a {\em weighted symmetric digraph}.
 A {\em dicycle} $C$ of $\G$ is a closed directed walk
 $(v_1,\ldots,v_{k+1})$ in which $v_1,\ldots,v_{k}$ are distinct
 vertices, $v_1=v_{k+1}$, and $v_iv_{i+1}$ ($i=1,\ldots,k$) are arcs.
 A {\em dipath} $P$ of $\G$ is a directed walk
 $(v_1,\ldots,v_{k})$ in which $v_1,\ldots,v_{k}$ are distinct
 vertices, and $v_iv_{i+1}$ ($i=1,\ldots,k-1$) are arcs.
 For a graph $G$ equipped with an orientation $\omega$, and
 a cycle $C$ of $G$ with a chosen direction of traversal (each cycle
 has two different directions for traversal), let $|C_\omega^+|$
 denote the number of edges of $C$ whose direction in $\omega$
 coincide with the direction of the traversal. Let $|C_\omega^-|$
 denote
 the value $|C|-|C_\omega^+|$ where $|C|$ is the length of $C$.
 Define $\tau(C,\omega)=\max\{|C|/|C_\omega^+|, |C|/|C_\omega^-|\}$.

 For reals $x$ and $r$, let $x($mod $r)$ be the unique value
 $t\in[0,r)$ such that $t\equiv x($mod $r)$. A {\em breaker function}
 of a weighted symmetric digraph $(\G,c)$ is a function $T:E(\G)\rightarrow \{0,1\}$ such that
 $T(xy)+T(yx)=1$ for each arc $xy\in E(\G)$. Henceforth,
 $T(xy)$ is written as $T_{xy}$.
 If $C$ is a dicycle of the weighted digraph $(\G,c)$ having a breaker
 function $T$, then the two values
 $\sum_{uv\in E(C)}T_{uv}$ and $\sum_{uv\in E(C)}c_{uv}$
 are denoted by $|C|_T$ and $|C|_c$ respectively, where $E(C)$ is
 the collection of all arcs in $C$. If $P$ is a dipath of $(\G,c)$
 then $|P|_T$ and $|P|_c$ are defined in the same way.

 Let $G$ be an undirected graph.
 A symmetric digraph $\G$ is said to be {\em derived
 from} $G$ if $V(G)=V(\G)$, and if $xy$ is an edge of $G$ then $xy$ and $yx$ are arcs of
 $\G$ and vice versa. If digraph $\G$ is derived from $G$ , then $G$
 is called the {\em underlying graph} of $\G$.

 Suppose $k\geq 2d\geq 1$ are positive integers.
 A $(k,d)$-$coloring$ of a
 graph $G$ is a mapping $f:V(G)\rightarrow
 \{0,1,\ldots,k-1\}$ such that for any edge $xy$ of $G$,
 $d\leq |f(x)-f(y)|\leq k-d$.
 The {\em circular chromatic number} $\chi_c(G)$ of $G$ is defined
 as
 \begin{center}
 $\chi_c(G)=\inf\{k/d:$ $G$ has a $(k,d)$-colorable $\}$.
 \end{center}

 For a real number $r\geq 1$, a {\em circular $r$-coloring} of a
 graph G is a function $f:V(G)\rightarrow [0,r)$ such that for any
 edge $xy$ of $G$, $1\leq |f(x)-f(y)|\leq r-1$. It is known
 \cite{survey, survey2005} that
 \begin{center}
 $\chi_c(G)=\inf\{r:$ $G$ has a circular $r$-coloring $\}$.
 \end{center}
 It is clear that $G$ has a $(k,d)$-coloring if and only if $G$ has
 a circular $k/d$-coloring.

 For a positive real $p$, let $S^p$ denote a circle with perimeter
 $p$ centered at the origin of ${\cal R}^2$. In the
 obvious way, we can identify the circle $S^p$ with the
 interval $[0, p)$. For $x,y\in S^p$,
 let $d_p(x,y)$ denote the length of the arc on $S^p$ from $x$ to $y$ in the
 clockwise direction if $x\not= y$, and let $d_p(x,y)=0$ if $x=y$.
 A {\em circular $p$-coloring} of a weighted digraph $(\vec{G},c)$
 is a function $\varphi: V(\G)\rightarrow S^p$ such that for each
 arc $uv$ of $\vec{G}$,
 $d_p(\varphi(u),\varphi(v))\geq c_{uv}$.
 The
 {\em circular chromatic number $\chi_c(\vec{G},c)$} of a
 weighted digraph $(\vec{G},c)$,
 recently introduced by Mohar \cite{mohar}, is defined as
 \begin{center}
 $\chi_c(\vec{G},c)=\inf\{p:$ $(\vec{G},c)$ has a circular
 $p$-coloring $\}$.
 \end{center}
 It is clear that $\chi_c(G)=\chi_c(\G,\1)$, where $\G$ is a digraph
 derived from $G$ and $\1(xy)=1$ for each arc $xy$ of $\G$.

 As we can see in \cite{yeh-zhu, survey, survey2005},
 the parameter $\chi_c(G)$ is a refinement of $\chi(G)$.
 The concept of circular chromatic number has attracted considerable
 attention in the past decade
 (see \cite{survey,survey2005} for a survey of research in this area).
 Further, the parameter
 $\chi_c(\vec{G},c)$ generalizes $\chi_c(G)$ from many application
 points of view. Readers are referred to \cite{yeh} for a connection of
 circular colorings of weighted digraphs and parallel computations.
 It was also shown in \cite{mohar} that
 the notion of $\chi_c(\vec{G},c)$
 generalizes
 the weighted circular colorings \cite{deuber},
 the linear arboricity of a graph and
 the metric traveling salesman problem.

 In this paper we consider weighted symmetric digraph, and explore the relation
 between circular coloring and breaker function of a weighted
 symmetric digraph.
 Our main result is Theorem \ref{main} which generalizes the work of Zhu
 \cite{orientation},
 a result of Tuza \cite{tuza}, and a result of Goddyn, Tarsi and Zhang \cite{gtz},
 all these results of
 \cite{orientation}, \cite{tuza}, and \cite{gtz} are
 generalizations of Minty's work in \cite{minty}.

 For a graph $G$ and a weighted symmetric digraph $(\G,c)$, let
 $M(G)$ (resp. $M(\G)$) denote the collection of all
 cycles (resp. dicycles) in $G$ (resp. $\G$), and let
 $L(\G,c)=\max$\{$c_{xy}+c_{yx}:xy$ is an arc of $\G$\}.
 In this paper, for our convenience, we say that a quotient has a value of
 infinity if its denominator is zero.


 \begin{theorem}
 \label{main}
 Let $(\G,c)$ be a weighted
 symmetric digraph. Suppose $r$ is a real number with $r\geq L(\G,c)$.
 Then $(\G,c)$ has a circular $r$-coloring if and only if
 $(\G,c)$ has a breaker
 function $T$ such that
 $$ \max_{C}{|C|_c\over |C|_T}\leq r, $$
 where the maximum is taken over all dicycle $C$ satisfying
 $0< |C|_c(\mbox{{\rm mod} $r$})< L(\G, c)$.
 \end{theorem}

 Theorem \ref{main} says that to show a weighted symmetric digraph $(\G,c)$ has a circular
 $r$-coloring, it suffices to check those dicycles $C$ of $\G$ for which $0< |C|_c(\mbox{mod $r$})<
 L(\G,c)$.
 Three simple consequences of Theorem \ref{main} are the following.


 \begin{corollary}
 \label{new}
 Suppose $r$ is a real number with $r\geq 2$.
 If $G$ has an orientation $\omega$ such that
 $ \max_{C}{\tau(C,\omega)}\leq r,$
 where the maximum is taken over all cycle $C$ satisfying
 $0< |C|(\mbox{mod $r$})< 2$,
 then $G$ has a circular $r$-coloring.
 \end{corollary}



 \begin{corollary}{\rm \cite{orientation}}
 \label{zhu}
 Suppose $k$ and $d$ are integers with $k\geq 2d\geq 1$.
 If $G$ has an orientation $\omega$ such that
 $ \max_{C}{\tau(C,\omega)}\leq k/d$,
 where the maximum is taken over all cycle $C$ satisfying
 $1\leq d|C|(\mbox{mod $k$})\leq 2d-1$,
 then $G$ has a $(k,d)$-coloring.
 \end{corollary}



 \begin{corollary}
 \label{new1}
 If for each cycle $C$ of $G$, $|C|(\mbox{\rm mod $r$})\not \in
 (0,2)$, then $G$ has a circular $r$-coloring.
 \end{corollary}

 It is clear that Tuza's result in \cite{tuza} is a special case of
 Corollary \ref{zhu} as $d=1$. It had been shown in \cite{orientation} that
 Minty's work of \cite{minty} and the following Goddyn, Tarsi and Zhang's
 result for $(k,d)$-coloring both are special cases of Corollary
 \ref{zhu}: if $G$ has an acyclic orientation $\omega$ such that
 $\max_{C\in M(G)}\tau(C,\omega)\leq k/d$, then $G$ has a $(k,d)$-coloring.
 In section \ref{main section} we prove the main result: Theorem
 \ref{main}, and derive the first two corollaries mentioned above.

 %
 %
 %
 %
 \section{The proof of the main result}
 \label{main section}


 \begin{proofif}
 Suppose
 $(\G,c)$ has a breaker
 function $T$ such that
 $$ \max_{C}{|C|_c\over |C|_T}\leq r,$$
 where the maximum is taken over all dicycle $C$ satisfying
 $0< |C|_c(\mbox{{\rm mod} $r$})< L(\G, c)$.
 Let $w$ be a function which assigns to each arc $xy$ of $\G$ a weight
 $w_{xy}=c_{xy}-rT_{xy}$.
 Let $G$ be the underlying graph of $\G$.
 For a spanning tree $\T$ of $G$
 and two vertices $x$ and $y$ in $G$, clearly there is a unique path $v_1v_2\ldots v_k$
 in $G$ having $v_1=x$, $v_k=y$, and $v_iv_{i+1}$ ($i=1,\ldots,k-1$)
 are edges of $\T$. The $x,y$-dipath $(v_1,\ldots,v_k)$ of $\G$ generated
 in this way is called {\em the dipath of $\G$ from $x$ to $y$ in
 $\T$}.
 Let $s$ be a fixed vertex in $G$.
 Given a spanning tree $\T$ of $G$, we define the function
 $f_\T:V(\G)\rightarrow \mathbf{R}$ as follows:
 \begin{itemize}
 \item $f_\T(s)=0$;
 \item If $x$ is a vertex other than $s$ then
 $f_\T(x)=\sum_{xy}w_{xy}$, where the summation is taken over all arcs in
 the dipath of $\G$ from $s$ to $x$ in $\T$.
 \end{itemize}
 The {\em weight of $\T$} is defined to be $\sum_{v\in V(\G)}f_\T(v)$ and
 is denoted by $f(\T)$.
 In the following, let $\T$ be a spanning tree of $G$ with maximum
 weight. Let $\varphi$ be a function which assigns to each vertex of
 $\G$ a color $f_\T(v)($mod $r$) in $[0,r)$ (and hence in $S^r$). We shall show that $\varphi$ is a
 circular $r$-coloring of $(\G,c)$. To prove this, let $xy$ and $yx$
 be a pair of arcs in $\G$ and consider the following cases.
 In these cases, we consider $\T$ as a rooted
 tree with root $s$, let $x'$ and $y'$ be the fathers of $x$ and $y$
 respectively.

 {\bf Case I.} Suppose that $x$ is not on the $s,y$-path of $\T$ and $y$ is
 not on the $s,x$-path of $\T$. Let $\T'$ be the spanning tree of $G$
 obtained from $\T$ by deleting its edge $x'x$ and adding the edge
 $xy$. Since $f(\T')\leq f(\T)$, we have $f_{\T'}(x)\leq f_{\T}(x)$,
 and hence $f_{\T}(y)+w_{yx}\leq f_{\T}(x)$ because $y$ is the
 father of $x$ in $\T'$. Then by symmetry we have  $f_{\T}(x)+w_{xy}\leq
 f_{\T}(y)$. Therefore $$w_{yx}\leq f_{\T}(x)-f_{\T}(y)\leq
 -w_{xy}.$$
 If $T_{xy}=1$ then we have $c_{yx}\leq f_{\T}(x)-f_{\T}(y)\leq r-c_{xy}$.
 If $T_{xy}=0$ then we have $c_{xy}\leq f_{\T}(y)-f_{\T}(x)\leq r-c_{yx}$.
 In both subcases we arrive at $d_r(\varphi(x),\varphi(y))\geq
 c_{xy}$, and hence $d_r(\varphi(y),\varphi(x))\geq
 c_{yx}$ by symmetry.

 {\bf Case II.} Suppose that either the $s,y$-path of $\T$ contains $x$
 or the $s,x$-path of $\T$ contains $y$. Since $xy$ and $yx$ both are arcs of $\G$,
 it suffices to
 consider the case that $y$ is on the $s,x$-path of $\T$. Let $P$ be
 the dipath of $\G$ from $y$ to $x$ in $\T$ and $C$ be the dicycle
 of $\G$ consisting of $P$ and the arc $xy$. Using the same method
 as in the previous case, we have $$w_{yx}\leq
 f_{\T}(x)-f_{\T}(y).$$ Note that $f_{\T}(x)=f_{\T}(y)+\sum_{uv\in
 E(P)}w_{uv}$. That is
 $$f_{\T}(x)-f_{\T}(y)=|P|_c-r|P|_T.$$

 Suppose $|P|_c\equiv \ell$ (mod $r$), where $0\leq \ell <r$.
 Consider the following three scenarios.

  {\bf Subcase II(a).} If $T_{xy}=1$ and $\ell >-w_{xy}$, then
  $r-c_{xy}<\ell <r$, and hence $0<\{|P|_c+c_{xy}\}$(mod $r$)$<c_{xy}$.
  It follows that $0<|C|_c$(mod $r$)$<c_{xy}$.
  We see at once that  $0< |C|_c(\mbox{{\rm mod} $r$})< L(\G, c)$.
  By the hypothesis, we have $|C|_c/|C|_T\leq r$
  which is equivalent to $|P|_c-r|P|_T\leq rT_{xy} -c_{xy}$.
  Therefore $f_{\T}(x)-f_{\T}(y)\leq -w_{xy}$.

  {\bf Subcase II(b).} If $T_{xy}=1$ and $\ell <w_{yx}$, then
  $0<\ell <c_{yx}$, and hence  $c_{xy}<|C|_c$(mod
  $r$)$<c_{xy}+c_{yx}$.
  We see at once that  $0< |C|_c(\mbox{{\rm mod} $r$})< L(\G, c)$.
  By the same argument as in the
  Subcase II(a), we arrive at $f_{\T}(x)-f_{\T}(y)\leq -w_{xy}$.

  {\bf Subcase II(c).} If $T_{xy}=1$ and $w_{yx}\leq \ell \leq
  -w_{xy}$, then $c_{yx}\leq \ell \leq r-c_{xy}$.
  Since $f_{\T}(x)-f_{\T}(y)=|P|_c-r|P|_T$ and $|P|_c=q'r+\ell$ for some integer $q'$,
  there exists an integer
  $q$ such that $f_{\T}(x)-f_{\T}(y)=qr+\ell$, and hence we
  have
  $d_r(\varphi(x),\varphi(y))\geq
  c_{xy}$ and $d_r(\varphi(y),\varphi(x))\geq
  c_{yx}$.

  Next suppose $-|P|_c\equiv \ell'$ (mod $r$), where $0\leq \ell' <r$.
 Consider the following three scenarios.

 {\bf Subcase II(d).} If $T_{xy}=0$ and $\ell' >-w_{yx}$, then
 $-r<-\ell'<c_{yx}-r$, and hence
 $-r+c_{xy}<-\ell'+c_{xy}<c_{xy}+c_{yx}-r$ which implies
 $c_{xy}<|C|_c$(mod $r$)$<c_{xy}+c_{yx}$. Now, by the same argument as in the
 Subcase II(b), we come to $f_{\T}(x)-f_{\T}(y)\leq -w_{xy}$.

 {\bf Subcase II(e).} If $T_{xy}=0$ and $\ell' <w_{xy}$, then
 $0<-\ell'+c_{xy}<c_{xy}$, and hence $0<\{|P|_c+c_{xy}\}$(mod $r$)$<c_{xy}$.
 By the same argument as in the
 Subcase II(a), we arrive at $f_{\T}(x)-f_{\T}(y)\leq -w_{xy}$.

 {\bf Subcase II(f).} If $T_{xy}=0$ and $w_{xy}\leq \ell' \leq
  -w_{yx}$, then $c_{yx}\leq r-\ell'\leq r-c_{xy}$. Similarly as in
  the Subcase II(c),
  since $f_{\T}(x)-f_{\T}(y)=|P|_c-r|P|_T$ and $-|P|_c=q'r+\ell'$ for some integer $q'$,
  there exists an integer
  $q$ such that $f_{\T}(x)-f_{\T}(y)=qr+(r-\ell')$.
  Now we can conclude that $d_r(\varphi(x),\varphi(y))\geq
  c_{xy}$ and $d_r(\varphi(y),\varphi(x))\geq
  c_{yx}$.

  As we can see from the above cases, no matter what is the value of
  $T_{xy}$, we always have $d_r(\varphi(x),\varphi(y))\geq
  c_{xy}$ and $d_r(\varphi(y),\varphi(x))\geq
  c_{yx}$. This completes the proof of the `if' part.
 \end{proofif}


 \begin{proofonlyif}
 Suppose that $(\G,c)$ has a circular $r$-coloring $\varphi:
 V(\G)\rightarrow [0,r)$. Note that here we view $[0,r)$ as $S^r$.
 We will show that $(\G,c)$ has a breaker
 function $T$ such that
 $$\max_{C}{|C|_c\over |C|_T}\leq r,$$
 where the maximum is taken over all dicycles $C$ of $\vec{G}$, which is a
 stronger result than what we state in Theorem \ref{main}.
 Define a mapping $T$ which assigns to each arc $xy$ of $\G$ a value
 from $\{0,1\}$ such that $T(xy)=1$ as $\varphi(x)>\varphi(y)$,
 and $T(xy)=0$ as $\varphi(x)<\varphi(y)$.
 Clearly, $T$ is a breaker function of $(\G,c)$.
 Consider the two possibilities, $T_{xy}=1$ and $T_{xy}=0$,
 separately. It is easy to check that for each arc $xy$ of $\G$ we
 have $\varphi(x)+c_{xy}\leq \varphi(y)+rT_{xy}$.

 Let $\hat{C}$ be a dicycle of $\G$ with
 ${|\hat{C}|_c/ |\hat{C}|_T}=\max_{C\in M(\G)}{|C|_c/
 |C|_T},$
 say $\hat{C}=(v_1,v_2,\ldots,v_k)$.
 From the result proved in the last paragraph, we have
 $$\varphi(v_i)+c_{v_iv_{i+1}}\leq \varphi(v_{i+1})+rT_{v_iv_{i+1}}\,\,\,(i=1,2,\ldots,k),$$
 where we let $v_{k+1}=v_1$. Adding up both side of the inequalities
 separately, we shall arrive at ${|\hat{C}|_c \leq r|\hat{C}|_T}$.
 Therefore $\max_{C\in M(\G)}{|C|_c/
 |C|_T}\leq r$, that completes the proof of the `only if' part.
 \end{proofonlyif}

 As a byproduct, the proof of Theorem \ref{main}
 also provides a direct proof of the
 following result which was proved by Mohar \cite{mohar} through using
 a result of Hoffman \cite{hoffman} and Ghouila-Houri \cite{gh}:
 $
 \chi_c(\vec{G},c)=\min_{T}\max_{C}{|C|_c/|C|_T},
 $
 where the minimum is taken over all breaker function $T$ of $(\vec{G},c)$
 and the maximum is taken over all dicycles $C$ of $\vec{G}$.

 \begin{proofcorollary2}
 Suppose $G$ has an orientation $\omega$ such that
 $ \max_{C}{\tau(C,\omega)}\leq r$,
 where the maximum is taken over all cycles $C$ of $G$
 satisfying $0< |C|(\mbox{mod $r$})< 2$.
 Let $\G$ be the symmetric digraph derived from $G$.
 Let $c$ denote the all $1$'s function on the arcs of $\G$.
 Define a mapping $T:E(\G)\rightarrow \{0,1\}$ in the following way:
 If $xy$ is an edge of $G$ oriented from $x$ to $y$ in the
 orientation $\omega$ then let $T_{xy}=1$ and $T_{yx}=0$.
 Clearly, $T$ is a
 breaker function of $(\G,c)$, and $L(\G, c)=2$.
 Thus the corollary follows from Theorem \ref{main} immediately.
 \end{proofcorollary2}

 \begin{proofcorollary3}
 Suppose $G$ has an orientation $\omega$ such that
 $ \max_{C}{\tau(C,\omega)}\leq k/d$,
 where the maximum is taken over all cycles $C$ of $G$
 satisfying $1\leq d|C|(\mbox{mod $k$})\leq 2d-1$.
 Let $r=k/d$. Let $C$ be a cycle of $G$ having $0< |C|(\mbox{mod
 $r$})<2$, that is $|C|=rq+\ell$ for some integer $q$ and real $\ell$
 satisfying
 $0<\ell<2$. Clearly, $d\ell$ is an integer such that $0<d\ell <2d\leq k$.
 It follows that $1\leq d|C|(\mbox{mod $k$})\leq 2d-1$. Our
 hypothesis implies that $\tau(C,\omega)\leq k/d$.
 Corollary \ref{new} now shows that $G$ has a circular
 $k/d$-coloring, and hence $G$ has a $(k,d)$-coloring
 (see Zhu's survey \cite{survey2005} for a proof of the last assertion).
 \end{proofcorollary3}

%
%

\end{document}